\documentclass[12pt]{article}

\usepackage{amsmath}
\usepackage{amsthm}
\usepackage{amssymb}
\usepackage{graphicx}
\usepackage{mathtools}
\usepackage{color}
\usepackage[dvipsnames]{xcolor}
\usepackage{hyperref}
\usepackage{algorithm}
\usepackage{algpseudocode}
\usepackage{comment}
\usepackage{caption}
\usepackage{subcaption}

\usepackage{setspace}
\usepackage[margin=1in]{geometry}

\DeclareMathOperator{\St}{St}
\DeclareMathOperator{\vecop}{vec}
\DeclareMathOperator{\matop}{mat}
\DeclareMathOperator{\tr}{trace}
\DeclareMathOperator{\st}{s.t.}
\DeclareMathOperator{\diag}{diag}

\newcommand{\R}{\mathbb{R}}
\newcommand{\Shor}{\text{\sc Shor}}
\newcommand{\BDiag}{\text{\sc DiagSum}}
\newcommand{\Kron}{\text{\sc Kron}}

\algrenewcommand\algorithmicrequire{\textbf{Input:}}
\algrenewcommand\algorithmicensure{\textbf{Output:}}

\author{Samuel Burer%
\thanks{Department of Business Analytics, University of Iowa,
Iowa City, IA, 52242-1994, USA. Email: {\tt samuel-burer@uiowa.edu}} \and
Kyungchan Park\thanks{Department of Business Analytics, University of Iowa,
Iowa City, IA, 52242-1994, USA. Email: {\tt kyungchan-park@uiowa.edu}.}%
}
\date{August 5, 2022}

\title{A Strengthened SDP Relaxation for Quadratic Optimization Over the
Stiefel Manifold}

\begin{document}

\maketitle

\begin{abstract}

\noindent We study semidefinite programming (SDP) relaxations for the
NP-hard problem of globally optimizing a quadratic function over the
Stiefel manifold. We introduce a strengthened relaxation based on two
recent ideas in the literature: (i) a tailored SDP for objectives with
a block-diagonal Hessian; (ii) and the use of the Kronecker matrix
product to construct SDP relaxations. Using synthetic instances on four
problem classes, we show that, in general, our relaxation significantly
strengthens existing relaxations, although at the expense of longer
solution times.

\end{abstract}

\section{Introduction}

Given positive integers $n$ and $p$, the {\em Stiefel manifold\/} is
the set of all real $n \times p$ matrices with orthonormal columns:
\[
\St(n,p):=\{U\in\mathbb{R}^{n\times p}:U^TU=I_p\},
\]
where $I_p$ denotes the $p\times p$ identity matrix. We study quadratic
optimization over $\St(n,p)$:
\begin{equation} \label{equ:main} \tag{QPS}
    \min \left\{ u^T H u + 2 \, g^T u : u = \vecop(U), \ U \in \St(n,p) \right\},
\end{equation}
where $u = \vecop(U)$ is the vector formed by stacking the columns
of $U$. Here, the data are the symmetric matrix $H \in \R^{np \times
np}$ and column vector $g \in \R^{np}$. We recover $U$ from $u$ via the
operator $U = \matop(u)$, which reassembles the subvectors of $u$.

(\ref{equ:main}) has many applications, including subspace tracking in
signal processing and control, the orthogonal Procrustes problem and
its generalizations, conjugate gradient for large eigenvalue problems,
clustering for data mining, and electronic structures computation
\cite{birtea2020second, chretien2020, edelman1998, elden1977, elden1999,
hu2020, kim2020, smith1993, wen2013, zhang2014}.

For $p=1$, (\ref{equ:main}) is the well-known trust-region subproblem,
which is polynomial-time solvable \cite{More.Sorensen.1983}. On
the other hand, problem (\ref{equ:main}) is NP-hard for general
$p$. In particular, \cite{nemirovski2007} (page 303) discusses
this computational complexity for the case $g \ne 0$ and $p = n$.
Accordingly, many researchers have focused on iterative methods
to obtain stationary points or local minimizers \cite{absil2008,
boumal2022intromanifolds} as well as convex relaxations of (\ref{equ:main}), e.g.,
via semidefinite programming, to obtain tight, valid lower bounds
on its optimal value \cite{anstreicher1999,
anstreicher2000, dodig2015, nemirovski2007, overton1992}. In this
paper, we are interested in deriving strengthened SDP relaxations of
(\ref{equ:main}).

As a quadratically constrained quadratic program, (\ref{equ:main}) has
a standard semidefinite programming (SDP) relaxation, called the {\em
Shor relaxation\/}, which we describe in Section \ref{sec:relaxations}
and denote as \Shor. In a recent paper, Gilman et al.~\cite{gilman2022}
studied instances of (\ref{equ:main}) with block-diagonal $H$, and
they adapted the Shor relaxation to this block structure. In addition,
they introduced a strengthening of the Shor relaxation, which proved
to be extremely effective for solving block-diagonal instances. In
Section \ref{sec:relaxations}, we further adapt their strengthening in a
straightforward manner to the case of general $H$, and we refer to this
relaxation as \BDiag.

Then we introduce an additional strengthening of \BDiag. We refer to
this strengthening as \Kron~since it is based on the {\em Kronecker-product
idea\/} introduced by Anstreicher \cite{anstreicher2017}. In particular,
Anstreicher showed how to derive a valid positive-semidefinite
constraint for the Shor relaxation of a quadratic program, which also
includes two linear matrix inequalities (LMIs), by considering the
Kronecker product of the LMIs. In our case, we will use the Kronecker
product of a single redundant LMI with itself. Full details are given in
Section \ref{sec:kron}.

We test the strength of \Kron~relative to \BDiag~and \Shor~on
synthetic instances for four problem classes, one of which consists of
block-diagonal $H$ similar to those considered in \cite{gilman2022}.
We show that \Kron~significantly strengthens \Shor~on all four problem
classes. It also improves \BDiag~on three of the four classes but does
not show an improvement on the block-diagonal instances. One downside
of \Kron~is that it requires significantly more time to solve than
\BDiag~and \Shor. We report results and solution times in Section
\ref{sec:NumExp}.

\subsection{Notation and terminology} \label{sec:notation}

Our notation and terminology are largely standard. We use the index $i$
to range over the set $\{1,\ldots,n\}$, and the indices $j$ and $k$ each
range over $\{1,\ldots,p\}$. In addition, we write
\[
    u_j := U_{\cdot j} \quad \forall \ j = 1,\ldots,p
    \quad\quad
    \Longrightarrow
    \quad\quad
    U = \begin{pmatrix} u_1 & \cdots & u_p \end{pmatrix},
    \quad
    u = \begin{pmatrix} u_1 \\ \vdots \\ u_p \end{pmatrix}.
\]
Note that the $(i,j)$-th component of $U$ equals the $i$-th component
of $u_j$, i.e., $U_{ij} = [u_j]_i$. For simplicity, we write $U_{ij} =
u_{ji}$, and we use pairs $(j,i)$ as a two-dimensional coordinate system
for accessing entries of vectors in $\mathbb{R}^{np}$ such as $u$. For
example, $u_{ji}$ refers to the entry of $u$ at position $(j - 1)n + i$.

\section{Three SDP Relaxations} \label{sec:relaxations}

In this section, we introduce the three SDP relaxations of
(\ref{equ:main}) mentioned in the Introduction: \Shor, \BDiag, and
\Kron. \Shor~is derived directly from (\ref{equ:main}) using standard
techniques from the semidefinite-programming literature. Then \BDiag~adds a valid constraint to \Shor, similar
to \cite{gilman2022}, and finally we construct \Kron~by adding a
new valid constraint to \BDiag. This new constraint is based on the
Kronecker-product idea introduced in \cite{anstreicher2017}. Hence,
the relaxations are ordered in terms of increasing (or more precisely,
non-decreasing) strength of their lower bounds on (\ref{equ:main}).
Our conceptual contribution in this paper is the observation that
the Kronecker-product idea applies in this setting, and we will
show empirically in Section \ref{sec:NumExp} that \Kron~significantly
improves the strength of \Shor~and \BDiag~in general.

\subsection{The Shor relaxation}

The Shor relaxation for (\ref{equ:main}), which we call \Shor, is
derived by relaxing the quadratic form
\[
    {1 \choose u}{1 \choose u}^T =
    \begin{pmatrix}
        1 \\ u_1 \\ \vdots \\ u_p
    \end{pmatrix}
    \begin{pmatrix}
        1 \\ u_1 \\ \vdots \\ u_p
    \end{pmatrix}^T
    \in \mathbb{S}^{1 + np}
\]
to the positive semidefinite matrix
\begin{equation} \label{equ:uXY_defn}
    Y := \begin{pmatrix} 1 & u^T \\ u & X \end{pmatrix} :=
    \begin{pmatrix}
        1 & u_1^T & \cdots & u_p^T \\
        u_1 & X_{11} & \cdots & X_{1p} \\
        \vdots & \vdots & \ddots & \vdots \\
        u_p & X_{p1} & \cdots & X_{pp}
    \end{pmatrix}
    \in \mathbb{S}^{1 + np}_+
\end{equation}
Note that each rank-1 outer product $u_j u_k^T \in \R^{n \times n}$
is relaxed to $X_{jk} \in \R^{n \times n}$. We will use a 0 index to
specify the first row and column of $Y$ as well as the $ji$-coordinate
system described in Section \ref{sec:notation} to access the remaining
rows and columns of $Y$. For example,
\[
    Y_{00} = 1, \quad\quad
    Y_{\cdot 0} = {1 \choose u}, \quad\quad
    Y_{\cdot ji} =
    \begin{pmatrix} u_{ji} \\ [X_{1j}]_{\cdot i} \\ \vdots \\ [X_{pj}]_{\cdot i} \end{pmatrix}
\]

In terms of $u = \vecop(U)$, the constraints defining $U \in \St(n,p)$
are $u_j^T u_j = 1$ for all $j$ and $u_j^T u_k = 0$ for all $j \ne k$.
Hence, \Shor~enforces $\tr(X_{jj}) = 1$ and $\tr(X_{jk}) = 0$. Moreover,
the objective function $u^T H u + 2 \, g^T u$ relaxes to $H \bullet X +
2 \, g^T u$, where $H \bullet X := \tr(HX)$, yielding
\begin{align}
    \min \ \ &H \bullet X + 2 \, g^T u \nonumber \\
    \st \ \ & \tr(X_{jj}) = 1 \quad \forall \ j = 1,\ldots,p \label{equ:shor} \tag{\Shor} \\
        &\tr(X_{jk}) = 0 \quad \forall \ j \ne k \nonumber \\
        &Y \succeq 0, \nonumber
\end{align}
where the optimization variables $u, X, Y$ are related according to the
definition (\ref{equ:uXY_defn}).

\subsection{The diagonal-sum relaxation} \label{sec:bdiag}

To state the relaxation \BDiag, which is based on \cite{gilman2022},
consider the implications
\begin{equation} \label{equ:schur}
    U^T U = I_p
    \quad \Longrightarrow \quad
    U^T U \preceq I_p
    \quad \Longleftrightarrow \quad
    \begin{pmatrix} I_p & U^T \\ U & I_n \end{pmatrix} \succeq 0
    \quad \Longleftrightarrow \quad
    UU^T \preceq I_n,
\end{equation}
where $I_n$ denotes the $n \times n$ identity matrix. The first
implication follows by relaxing the Stiefel equation $U^T U = I_p$ to a
matrix inequality, and the second and third implications follow by the
Schur complement theorem. Then
\[
    \sum_{j=1}^p u_j u_j^T = UU^T \preceq I_n
    \quad \Longrightarrow \quad
    \sum_{j=1}^p X_{jj} \preceq I_n,
\]
and \BDiag~is formed by adding this final linear matrix inequality to
the Shor relaxation:
\begin{align}
    \min \ \ &H \bullet X + 2 \, g^T u \nonumber \\
    \st \ \ & \tr(X_{jj}) = 1 \quad \forall \ j = 1,\ldots,p \nonumber \\
        &\tr(X_{jk}) = 0 \quad \forall \ j \ne k \label{equ:bdiag} \tag{\BDiag} \\
        &Y \succeq 0 \nonumber \\
        &\sum_{j=1}^p X_{jj} \preceq I_n. \nonumber
\end{align}

In fact, \cite{gilman2022} studied instances of (\ref{equ:main}) for
which $H$ is a block-diagonal matrix of the form $H := \text{Diag}\left(
H_{11}, \ldots, H_{pp} \right)$ with each $H_{jj} \in \mathbb{S}^n$. In
this case, the objective function simplifies to $\sum_{j=1}^p H_{jj}
\bullet X_{jj} + 2 \, g^T u$. The precise form of the relaxation
considered in \cite{gilman2022} then dropped the off-diagonal $X_{jk}$,
considered the weaker semidefinite conditions
\[
    \begin{pmatrix} 1 & u_j^T \\ u_j & X_{jj} \end{pmatrix} \succeq 0 \quad
    \forall \ j = 1, \ldots, p
\]
corresponding to certain principal submatrices of $Y \succeq 0$, and
enforced $\sum_{j=1}^p X_{jj} \preceq I_n$.\footnote{Note that, although
the relaxation in \cite{gilman2022} was applied to $g \ne 0$, the authors
only experimented with $g = 0$ due to their specific application of
interest. We will also test with $g = 0$ in Section \ref{sec:NumExp}.}

\subsection{The Kronecker relaxation} \label{sec:kron}

To derive \Kron, we observe that the Kronecker-product idea developed in
\cite{anstreicher2017} can be applied in this setting. Specifically,
using (\ref{equ:schur}) and the fact that the Kronecker product of
positive semidefinite matrices is positive semidefinite, we see that
\begin{equation} \label{equ:kron}
    \begin{pmatrix} I_p & U^T \\ U & I_n \end{pmatrix} \succeq 0
    \quad \Longrightarrow \quad
    \begin{pmatrix} I_p & U^T \\ U & I_n \end{pmatrix}
    \otimes
    \begin{pmatrix} I_p & U^T \\ U & I_n \end{pmatrix}
    \succeq 0.
\end{equation}
The resultant linear matrix inequality depends quadratically on
$U$, and hence it can be linearized in the space of the lifted
variable $Y$. Note that the size of the Kronecker-product matrix is
$(n + p)^2 \times (n + p)^2$.

To determine the precise form of the linearization, define the constant
matrices
\[
    K_{ji} = e_{p+i} e_j^T + e_j e_{p+i}^T \in \mathbb{S}^{p + n}
    \quad\quad \forall \ \ (j,i) \in \{1,\ldots,p\} \times \{1,\ldots,n\}.
\]
We then have the expression
\begin{align*}
    \begin{pmatrix}
        I_p & U^T \\ U & I_n
    \end{pmatrix}
    &= I_{p+n} + \sum_{i=1}^n \sum_{j=1}^p U_{ij} (e_{p+i} e_j^T + e_j e_{p+i}^T) \\
    &= I_{p+n} + \sum_{j=1}^p \sum_{i=1}^n u_{ji} K_{ji}.
\end{align*}
Hence,
\begin{align*}
    \begin{pmatrix} I_p & U^T \\ U & I_n \end{pmatrix}
    &\otimes
    \begin{pmatrix} I_p & U^T \\ U & I_n \end{pmatrix} =
    I_{p+n} \otimes I_{p+n} + \\
    &\sum_{j=1}^p \sum_{i=1}^n u_{ji} ( I_{p + n} \otimes K_{ji} + K_{ji} \otimes I_{p+n}) +
    \sum_{j=1}^p \sum_{i=1}^n \sum_{k=1}^p \sum_{l=1}^n
    u_{ji} u_{kl} K_{ji} \otimes K_{kl},
\end{align*}
and so the quadratic linear matrix inequality introduced above
linearizes to $M(u,X) \succeq 0$, where
\begin{align*}
    M(u, X) := I_{(p+n)^2} &+
    \sum_{j=1}^p \sum_{i=1}^n u_{ji} ( I_{p + n} \otimes K_{ji} + K_{ji} \otimes I_{p+n}) + \\
                           &\sum_{j=1}^p \sum_{i=1}^n \sum_{k=1}^p \sum_{l=1}^n
    [X_{jk}]_{il} K_{ji} \otimes K_{kl}.
\end{align*}

Adding $M(u,X) \succeq 0$ to \BDiag, we obtain the \Kron~relaxation
\begin{align}
    \min \ \ &H \bullet X + 2 \, g^T u \nonumber \\
    \st \ \ & \tr(X_{jj}) = 1 \quad \forall \ j = 1,\ldots,p \nonumber \\
        &\tr(X_{jk}) = 0 \quad \forall \ j \ne k \label{equ:bdiag} \tag{\Kron} \\
        &Y \succeq 0 \nonumber \\
        &\sum_{j=1}^p X_{jj} \preceq I_n \nonumber \\
        &M(u,X) \succeq 0. \nonumber
\end{align}

\subsubsection{Two variations}

Before taking the Kronecker product in (\ref{equ:kron}), we could also
symmetrically permute one of the two matrices to
\[
    \begin{pmatrix} I_n & U \\ U^T & I_p \end{pmatrix} \succeq 0
\]
and consider the product
\[
    \begin{pmatrix} I_n & U \\ U^T & I_p \end{pmatrix}
    \otimes
    \begin{pmatrix} I_p & U^T \\ U & I_n \end{pmatrix}
    \succeq 0.
\]
However, after linearization, it is not difficult to see that this would
simply lead to a permutation of $M(u,X) \succeq 0$ and hence would not
generate a new valid constraint. In fact, any symmetric permutation will
lead to the same constraint $M(u,X) \succeq 0$.

Instead of considering the Kronecker product, we could use the fact
that the Hadamard product (indicated by the ``$\circ$'' symbol) of positive semidefinite matrices is positive
semidefinite to relax the $(n + p) \times (n + p)$ condition
\[
    \begin{pmatrix} I_p & U^T \\ U & I_n \end{pmatrix}
    \circ
    \begin{pmatrix} I_p & U^T \\ U & I_n \end{pmatrix}
    \succeq 0
\]
to
\begin{equation} \label{equ:local1}
    \begin{pmatrix}
        I_p & \matop(\diag(X))^T \\
        \matop(\diag(X)) & I_n
    \end{pmatrix} \succeq 0.
\end{equation}
The paper \cite{jiang2019} considered such an approach in a more general
setting. It is clear that this positive semidefinite matrix is a small
principal submatrix of $M(u,X) \succeq 0$, and so (\ref{equ:local1})
could be added to \BDiag~to form a fourth relaxation, whose strength
is in between \BDiag~and \Kron~and could be solved more quickly
than \Kron. Through extensive numerical experiments, however, we could
find no instances for which (\ref{equ:local1}) strengthens \BDiag.
Said differently, it appears that \BDiag~implies (\ref{equ:local1}),
although we have so far not been able to prove this. Due to this
experience, we will not consider the Hadamard product further in this
paper.

\subsection{Dual bounds and primal values} \label{sec:dpvals}

For each of the three SDP relaxations developed in the preceding
subsections, its optimal value $d$ is a lower bound on the optimal value
$v$ of (\ref{equ:main}), i.e., $d \le v$. We can also use any feasible
solution $(u,X)$ of the SDP to generate a feasible value (or primal
value) $p$ for (\ref{equ:main}) such that $d \le v \le p$. This provides
an absolute duality gap $p - d$ on the optimal value $v$, or a relative
gap $(p - d)/\max\{1, |\tfrac{1}{2} (p + d)|\}$.

Indeed, given feasible $(u,X)$ and $Y$ related by (\ref{equ:uXY_defn}),
if $\text{rank}(Y)$ happens to equal 1, then by construction $U
:= \matop(u)$ is an element of $\St(n,p)$, and then $p := u^T H u
+ 2 \, g^T u$ is the desired primal value. On the other hand, if
$\text{rank}(Y) > 1$, we can still calculate a primal value as follows:

\begin{itemize}

    \item let $U_0 \Sigma_0 V_0^T$ be the thin singular value
    decomposition of $U$ and define $\tilde U := U_0 V_0^T \in \St(n,p)$;

    \item then set $\tilde u := \vecop(\tilde U)$ and $p := \tilde u^T H
    \tilde u + 2 \, g^T \tilde u$.

\end{itemize}

\noindent In the numerical experiments of Section \ref{sec:NumExp}, we
calculate $p$ from the optimal solution $(u,X)$ of the SDP. In addition,
we will use the solver Manopt \cite{manopt}, which is a toolbox for
local optimization over manifolds, to improve heuristically the value
$p$ starting from the matrix $\tilde U$. Note that Manopt contains
built-in data structures to handle the Stiefel manifold constraint.
In particular, we run Manopt's {\tt trustregions} algorithm over
the Stiefel manifold with its default settings. We need only supply
function, gradient, and Hessian evaluations for our objective function,
which are straightforward given that the function is quadratic.

\section{Numerical Experiments} \label{sec:NumExp}

In this section, we evaluate the quality of the dual bounds and primal
values provided by \Kron~in comparison with \BDiag~and~\Shor. We also
investigate the time required to solve the various SDP relaxations. We
test various combinations of $(n,p)$ ranging from $(6,2)$ to $(12,11)$.

In preliminary tests, we also evaluated the global solver Gurobi
version 9.5.1, which solves (\ref{equ:main}) by calculating a sequence
of primal values and dual bounds with gaps converging to zero. With
$(n,p)=(6,3)$, for example, the gaps calculated by Gurobi after 15
minutes were significantly worse than those calculated by \Shor,
\BDiag, and \Kron~each within 1 second. Moreover, Gurobi's time and gap
performance for larger values of $(n,p)$ was even worse relative to
the three SDP relaxations. Hence, in what follows, we do not include
explicit comparisons with Gurobi.

We first describe our problem instances and then specify our testing
algorithm and measurements. Finally, we summarize our results.

\subsection{Problem instances} \label{sec:probs}

% We remark that, throughout the results shown below in Section
% \ref{sec:results}, we will take $n = 10$ and $p = 5$ because these
% values are good representatives of the general case. On the other
% hand, we now describe problem classes, where $n$ and $p$ are given
% arbitrarily.

In this subsection, we consider positive integers $n \ge p$ to be fixed.
We consider four problem classes:

\begin{enumerate}

\item Our first class of instances is based on random $H \in
\mathbb{S}^{np}$ and $g \in \mathbb{R}^{np}$ such that every entry in
$H$ and $g$ is i.i.d.~${\cal N}(0,1)$. We call these the {\em random
instances\/}.

\item Our second class of instances is the {\em block-diagonal
instances\/} in which $H$ is a block-diagonal matrix of the form $H :=
\text{Diag}\left( H_{11}, \ldots, H_{pp} \right)$ with each $H_{jj}
\in \mathbb{S}^n$. We also take $g = 0$ so that $H \bullet X + 2 \,
g^T u = \sum_{j=1}^p H_{jj} \bullet X_{jj}$. As discussed in Section
\ref{sec:bdiag}, this class was the basis of the numerical experiments
in \cite{gilman2022} for a subclass of matrices $\{ H_{jj} \}$
satisfying certain structural properties. For example, their matrices
$H_{jj}$ were negative semidefinite. The authors found that \BDiag~was
typically very strong for these instances, but some instances still
showed a positive gap. In our case, we randomly generate all entries of
$\{ H_{jj} \}$ by drawing each entry i.i.d.~${\cal N}(0,1)$.

\item Our third class is the {\em Procrustes instances\/},
which are random instances of the orthogonal Procrustes problem
\cite{gower2004procrustes}. Given an additional positive integer $m$ and
matrices $A\in\mathbb{R}^{m \times n}$ and $B\in\mathbb{R}^{m \times
p}$, the problem is to find an orthogonal mapping of the columns of $A$
into $\mathbb{R}^{m \times p}$, which minimizes the Frobenius distance
to $B$, i.e.,
\[
\displaystyle\min_{U\in\St(n,p)} \lVert AU-B\rVert_F^2,
\]
where the Frobenius norm is defined by $\|M\|_F := \sqrt{M \bullet M}$.
This is an instance of (\ref{equ:main}) with
\[
    H := I_p\otimes(A^T A),  \quad g := \vecop(-A^T B).
\]
Based on $n$ and $p$, to generate a Procrustes instance, we choose
$m$ to be uniformly distributed between $\lceil n/2 \rceil$ and
$2n$, inclusive, and then we randomly generate $A$ and $B$ with
entries i.i.d.~${\cal N}(0,1)$. We remark that there exist several
special cases, which are known to have closed-form solutions; see
\cite{elden1999,breloy2021}. Three of these cases are: $A = I_n$ with $m
= n$; $n = p$; and $p = 1$.

\item Our fourth and final class is the {\em Penrose instances\/},
which are random instances of the Penrose regression problem, which
is itself a generalization of the orthogonal Procrustes problem; see
\cite{chu1997orthogonally, elden1999, birtea2020second, kim2020}.
Given two additional positive integers $m$ and $q$, and matrices
$A\in\mathbb{R}^{m\times n}$, $B\in\mathbb{R}^{m\times q}$ and
$C\in\mathbb{R}^{p\times q}$, the problem is
\[
\min_{U\in\St(n,p)} \lVert AUC-B\rVert_F^2,
\]
which is an instance of (\ref{equ:main}) with
\[
    H := (CC^T)\otimes(A^T A), \quad
    g := vec(-A^T BC^T).
\]
Based on $n$ and $p$, we choose both $m$ and $q$ to be uniformly and
independently distributed between $ \lceil n/2 \rceil$ and $2n$,
inclusive, and then we randomly generate $A$, $B$, and $C$ with entries
i.i.d.~${\cal N}(0,1)$. Note that, by setting $C=I_p$ with $q=p$, we
recover the Procrustes problem as a special case.

\end{enumerate}

\subsection{Testing algorithm, measurements, and experiments}

Given a single instance of (\ref{equ:main}), Algorithm \ref{alg:sdp}
specifies our testing procedure on that instance. In words, we
solve $\Shor$, $\BDiag$, and $\Kron$~on the instance and save the
corresponding relative gaps and solution times.

\begin{algorithm}
    \caption{Approximate a Single Instance of (\ref{equ:main})}\label{alg:sdp}
    \begin{algorithmic}[1] %\mathrlap{}\hphantom{}  \hphantom{}\mathllap{}
        \Require An instance of (\ref{equ:main}), i.e., positive
        integers $n \ge p$ and data $H \in \mathbb{S}^{np}$ and $g \in \mathbb{R}^{np}$.
        % \State Calculate (heuristic) $p_m$ based on $np$ random starting points.
        \For{$\text{{\sc rlx}} \in \{ \Shor, \BDiag, \Kron \}$}
            \State Solve {\sc rlx} to obtain an optimal solution
            $(u^*, X^*)$.
            \State Save the total time $t_{\text{{\sc rlx}}}$ (in
            seconds) for optimizing {\sc rlx}.
            \State Set $d_{\text{{\sc rlx}}} = H \bullet X^* + 2 \, g^T u^*$.
            \State Calculate $p_{\text{{\sc rlx}}}$ based on $(u^*,X^*)$ as
            described in Section \ref{sec:dpvals}.
            \State Calculate the relative gap $\gamma_{\text{{\sc
            rlx}}} := (p_{\text{{\sc rlx}}} - d_{\text{{\sc rlx}}}) /
            \max\{1,
                |\tfrac12 (p_{\text{{\sc rlx}}} + d_{\text{{\sc rlx}}})|
            \}$.
        \EndFor
        \Ensure Relative gaps $\gamma_{\Shor}, \gamma_{\BDiag},
        \gamma_{\Kron}$ and times $t_{\Shor}, t_{\BDiag},
        t_{\Kron}$.
    \end{algorithmic}
\end{algorithm}

\noindent Based on the output of Algorithm \ref{alg:sdp}, we say
that a given relaxation {\sc rlx} {\em solves\/} the instance if
the relative gap $\gamma_{\text{\sc rlx}}$ is less than $10^{-4}$.

For each of the three values $n \in \{6, 9, 12\}$, we tested the three values
$p \in \{2, \lceil n/2 \rceil, n-1\}$ for a total of nine pairs $(n,p)$
ranging from $(6,2)$ to $(12,11)$. For a fixed pair $(n,p$), we ran
Algorithm \ref{alg:sdp} on 1,000 instances of each of the four problem
classes described in Section \ref{sec:probs}.

All experiments were coded in MATLAB (version 9.11.0.1873467, R2021b
Update 3), and all SDP relaxations were solved with Mosek (version
9.3.10). Manopt version 7.0 was employed for the local optimization
subroutine described in Section \ref{sec:dpvals}. The computing
environment was a single Xeon E5-2680v4 core running at 2.4 GHz with 8
GB of memory under the CentOS Linux operating system.

For \Shor~and \BDiag, we call Mosek directly, whereas we construct
\Kron~using the modeling interface YALMIP (version 31-March-2021), which
then calls Mosek. YALMIP is used for coding simplicity, since the
Kronecker constraint was relatively challenging to implement directly
in Mosek. Furthermore, we verified in two ways that the model passed
to YALMIP from Mosek was an efficient representation of \Kron. First,
we used YALMIP's {\tt dualize} command to test the dual formulation;
it required significantly more time to solve. Second, we verified that
the problem size reported by Mosek matched what one would expect from a
direct implementation of \Kron.

\subsection{Results} \label{sec:results}

We first summarize the solution times in Table \ref{tab:times}, which
shows the average times (in seconds) for the three SDP relaxations
over all instances, grouped by the nine $(n,p)$ pairs. It is clear that
\Kron~requires significantly more time than both \Shor~and \BDiag~as $n$
and $p$ increase. \BDiag~also requires more time than \Shor, but its
growth with $n$ and $p$ is much more modest than \Kron's.

\begin{table}[ht]
\centering
\begin{tabular}{rrccr}
$n$ & $p$ & \Shor & \BDiag & \Kron \\ \hline\hline
6 &    2 & 0.011 & 0.031 & 0.626 \\
6 &    3 & 0.016 & 0.040 & 1.027 \\
6 &    5 & 0.037 & 0.108 & 3.008 \\ \hline
9 &    2 & 0.015 & 0.066 & 2.284 \\
9 &    5 & 0.198 & 0.359 & 11.645 \\
9 &    8 & 0.618 & 0.864 & 47.000 \\ \hline
12 &    2 & 0.025 & 0.142 & 6.297 \\
12 &    6 & 0.584 & 0.875 & 54.389 \\
12 &   11 & 0.707 & 0.813 & 250.233 \\
\end{tabular}
\caption{Average solution times (in seconds) for the three SDP
relaxations over all instances, grouped by the nine $(n,p)$ pairs}
\label{tab:times}
\end{table}

\begin{figure}
    \centering
    \includegraphics[width=6.5in]{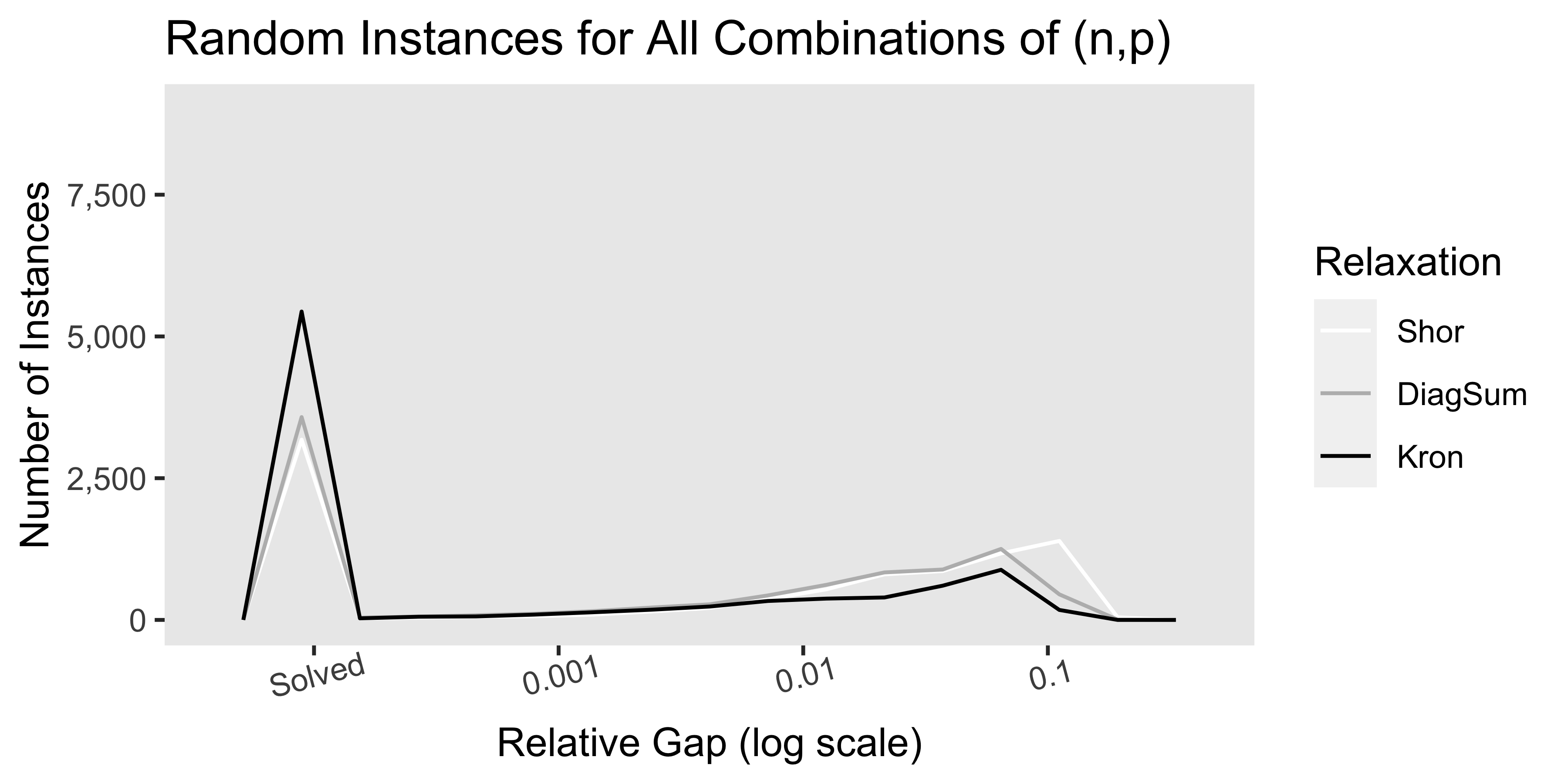}
    \caption{Histograms (depicted as line charts) of the relative gaps
    for the three SDP relaxations on 9,000 random instances}
   \label{fig:random}
\end{figure}

\begin{figure}
    \centering
    \includegraphics[width=6.5in]{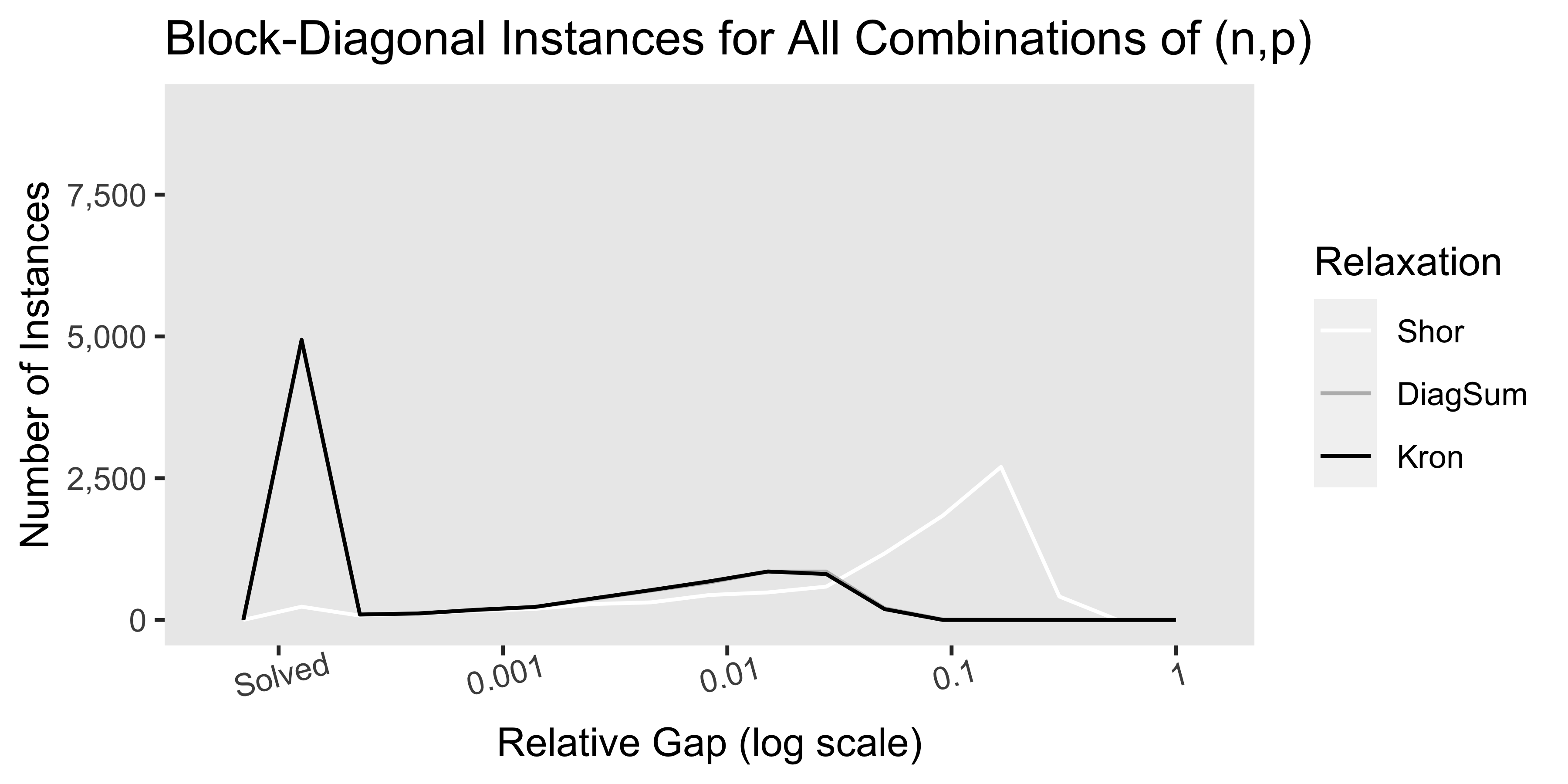}
    \caption{Histograms (depicted as line charts) of the relative gaps
    for the three SDP relaxations on 9,000 block-diagonal instances}
   \label{fig:block-diagonal}
\end{figure}

\begin{figure}
    \centering
    \includegraphics[width=6.5in]{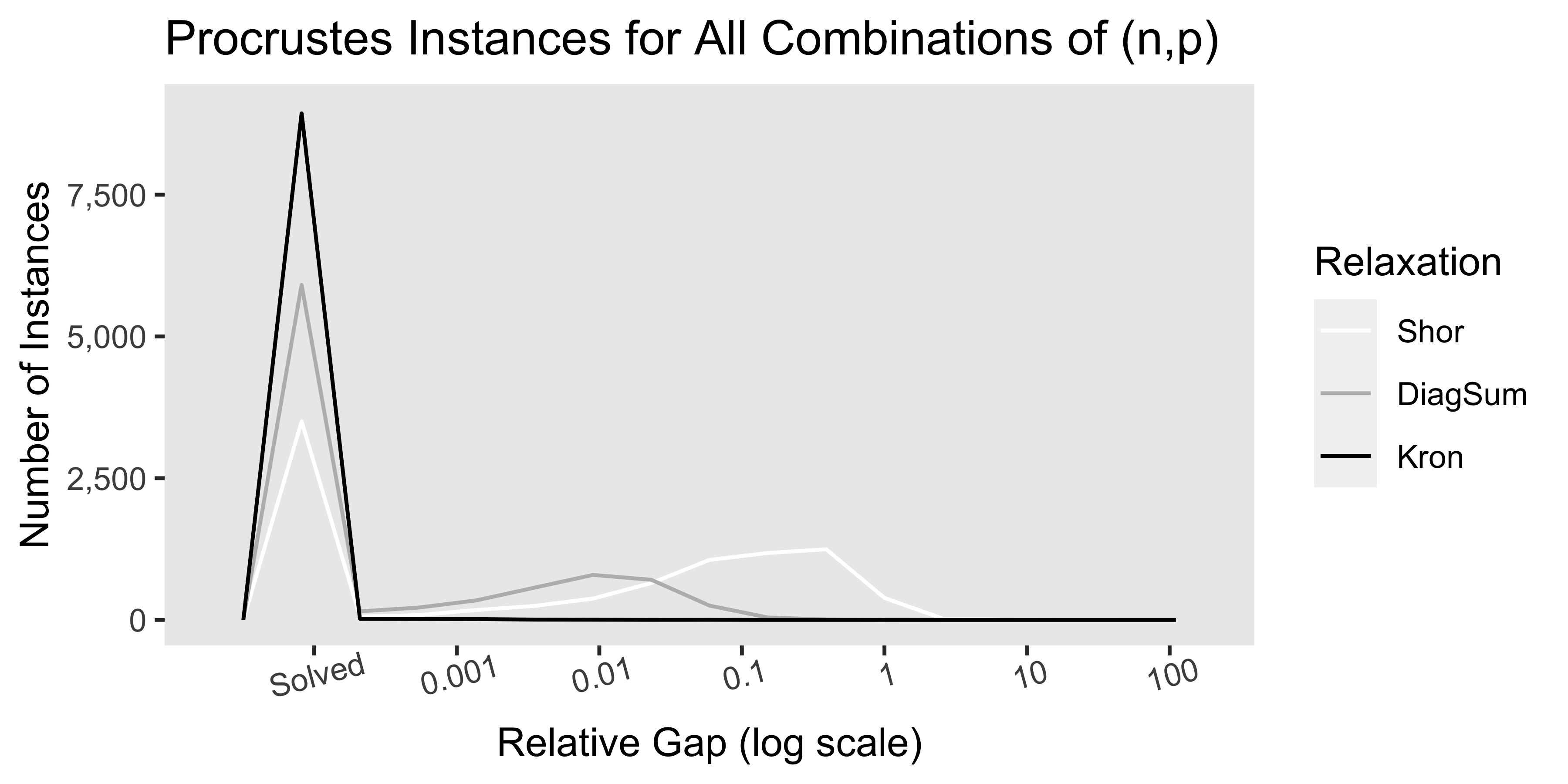}
    \caption{Histograms (depicted as line charts) of the relative gaps
    for the three SDP relaxations on 9,000 Procrustes instances}
   \label{fig:procrustes}
\end{figure}

\begin{figure}
    \centering
    \includegraphics[width=6.5in]{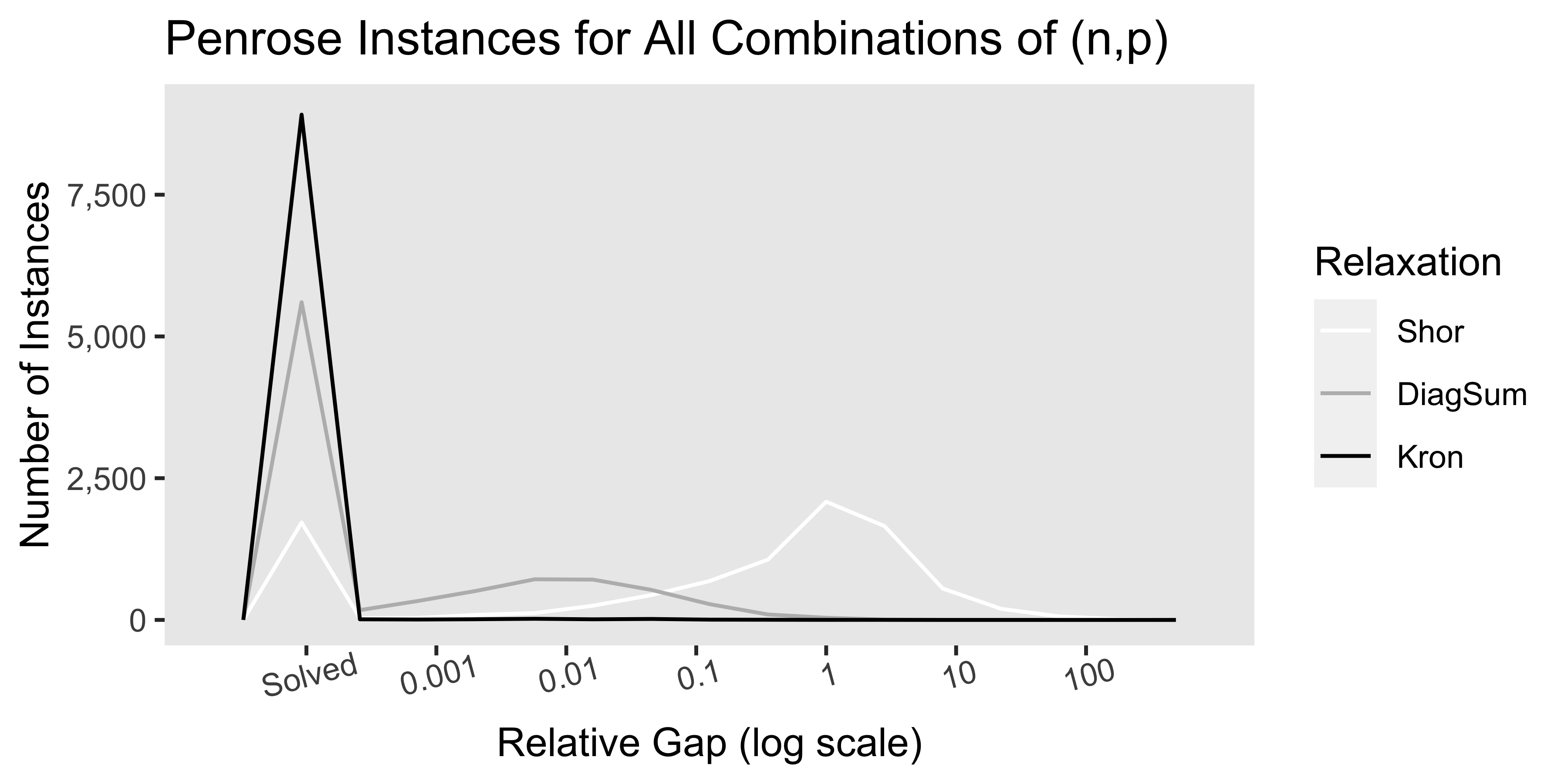}
    \caption{Histograms (depicted as line charts) of the relative gaps
    for the three SDP relaxations on 9,000 Penrose instances}
   \label{fig:penrose}
\end{figure}

Next, we examine the relative gaps in Figures
\ref{fig:random}--\ref{fig:penrose}. Each figure corresponds
to one problem class and
shows the histograms of the relative gaps (depicted as line
charts, or ``frequency polygons,'' for visibility) for \Shor, \BDiag~and
\Kron~over all instances of that problem class---9,000 instances in
total since there are 1,000 instances for each of the nine $(n,p)$
pairs. As discussed in Section \ref{sec:dpvals}, an instance is {\em
solved\/} by a relaxation if the relative gap is less than $10^{-4}$.
Such instances are grouped in the {\em Solved\/} category in the figure.

% \begin{figure}
%     \centering
%     \includegraphics[width=6.5in]{fig_all}
%     \caption{Histograms (depicted as line charts) of the relative gaps
%     for the three SDP relaxations on all 24,000 instances}
%    \label{fig:all}
% \end{figure}

% Figures \ref{fig:random}--\ref{fig:penrose} depict the same information
% except segmented by problem class, each figure representing 9,000
% instances in total. The random instances in Figure \ref{fig:random}
% depict the same overall pattern seen in Figure \ref{fig:all}.

For the random instances in Figure \ref{fig:random}, it is clear
that \Kron~is the strongest relaxation with more instances solved
and generally lower relative gaps. On the other hand, the remaining
figures show different patterns. For the block-diagonal instances
in Figure \ref{fig:block-diagonal}, \Kron~showed little to no
improvement over \BDiag, whereas both relaxations significantly
strengthened \Shor. For the Procrustes and Penrose instances in Figures
\ref{fig:procrustes}--\ref{fig:penrose}, the performance of \Kron~is
dramatically better than both \Shor~and \BDiag, where \Kron~solved
nearly all instances. (Although it is difficult to see in the figures,
there were a few instances {\em not\/} solved by \Kron.)

\section{Conclusions}

We have shown how the Kronecker-product idea can be
applied to improve SDP relaxations of (\ref{equ:main}), although the
key constraint---the linear matrix inequality $M(u,X) \succeq 0$---adds
to the complexity and time required to solve
the relaxation. An area of future research would be to incorporate
the strength of the Kronecker constraint without incurring the full
computational cost, e.g., by a clever cutting plane method or by
enforcing positive semidefiniteness on critical principal submatrices of
$M(u,X)$.

Given the small relative gaps in Figure
\ref{fig:procrustes}--\ref{fig:penrose} on the Procrustes and Penrose
instances, it is worthwhile to investigate whether one can prove a
guaranteed quality of \Kron~on these problem classes. Numerically,
we found that \Kron~did {\em not\/} solve all instances exactly, but
given that numerical computation is necessarily inexact and relies
on, for example, the accuracy of the underlying SDP solver, it would
be interesting to either exhibit analytically an instance for which
\Kron~is inexact or to prove that \Kron~is tight for the Procrustes or
Penrose classes.

Another avenue for improvement would be to extend the theoretical
results of Gilman et al \cite{gilman2022}, from which we adapted the
\BDiag~relaxation. Indeed, their work provides a global optimality
certificate for a block-diagonal instance of (\ref{equ:main}) by solving
a dimension-reduced SDP feasibility problem, which is closely related to
their SDP relaxation. Perhaps these results could be brought to bear on
our problem for general $H$.

\section*{Acknowledgements}

We thank Kurt Anstreicher and Laura Balzano for contributing extremely
helpful comments on a draft version of this research.

\bibliographystyle{abbrv}
\bibliography{paper}

\begin{thebibliography}{10}

\bibitem{absil2008}
P.-A. Absil, R.~Mahony, and R.~Sepulchre.
\newblock {\em Optimization algorithms on matrix manifolds}.
\newblock Princeton University Press, Princeton, NJ, 2008.
\newblock With a foreword by Paul Van Dooren.

\bibitem{anstreicher1999}
K.~Anstreicher, X.~Chen, H.~Wolkowicz, and Y.-X. Yuan.
\newblock Strong duality for a trust-region type relaxation of the quadratic
  assignment problem.
\newblock {\em Linear Algebra Appl.}, 301(1-3):121--136, 1999.

\bibitem{anstreicher2000}
K.~Anstreicher and H.~Wolkowicz.
\newblock On lagrangian relaxation of quadratic matrix constraints.
\newblock {\em SIAM Journal on Matrix Analysis and Applications}, 22(1):41--55,
  2000.

\bibitem{anstreicher2017}
K.~M. Anstreicher.
\newblock Kronecker product constraints with an application to the
  two-trust-region subproblem.
\newblock {\em SIAM Journal on Optimization}, 27(1):368--378, 2017.

\bibitem{birtea2020second}
P.~Birtea, I.~Ca{\c{s}}u, and D.~Com{\u{a}}nescu.
\newblock Second order optimality on orthogonal stiefel manifolds.
\newblock {\em Bulletin des Sciences Math{\'e}matiques}, 161:102868, 2020.

\bibitem{boumal2022intromanifolds}
N.~Boumal.
\newblock An introduction to optimization on smooth manifolds.
\newblock To appear with Cambridge University Press, Apr 2022.

\bibitem{manopt}
N.~Boumal, B.~Mishra, P.-A. Absil, and R.~Sepulchre.
\newblock {M}anopt, a {M}atlab toolbox for optimization on manifolds.
\newblock {\em Journal of Machine Learning Research}, 15(42):1455--1459, 2014.

\bibitem{breloy2021}
A.~Breloy, S.~Kumar, Y.~Sun, and D.~P. Palomar.
\newblock Majorization-minimization on the stiefel manifold with application to
  robust sparse pca.
\newblock {\em IEEE Transactions on Signal Processing}, 69:1507--1520, 2021.

\bibitem{chretien2020}
S.~Chr{\'e}tien and B.~Guedj.
\newblock Revisiting clustering as matrix factorisation on the stiefel
  manifold.
\newblock In {\em International Conference on Machine Learning, Optimization,
  and Data Science}, pages 1--12. Springer, 2020.

\bibitem{chu1997orthogonally}
M.~T. Chu and T.~T. Nickolay.
\newblock The orthogonally constrained regression revisited.
\newblock 1997.

\bibitem{dodig2015}
M.~Dodig, M.~Sto{\v{s}}i{\'c}, and J.~Xavier.
\newblock On minimizing a quadratic function on stiefel manifold.
\newblock {\em Linear Algebra and its Applications}, 475:251--264, 2015.

\bibitem{edelman1998}
A.~Edelman, T.~A. Arias, and S.~T. Smith.
\newblock The geometry of algorithms with orthogonality constraints.
\newblock {\em SIAM journal on Matrix Analysis and Applications},
  20(2):303--353, 1998.

\bibitem{elden1977}
L.~Eld{\'e}n.
\newblock Algorithms for the regularization of ill-conditioned least squares
  problems.
\newblock {\em BIT Numerical Mathematics}, 17(2):134--145, 1977.

\bibitem{elden1999}
L.~Eld{\'e}n and H.~Park.
\newblock A procrustes problem on the stiefel manifold.
\newblock {\em Numerische Mathematik}, 82(4):599--619, 1999.

\bibitem{gilman2022}
K.~Gilman, S.~Burer, and L.~Balzano.
\newblock A semidefinite relaxation for sums of heterogeneous quadratics on the
  stiefel manifold, 2022.

\bibitem{gower2004procrustes}
J.~C. Gower and G.~B. Dijksterhuis.
\newblock {\em Procrustes problems}, volume~30.
\newblock OUP Oxford, 2004.

\bibitem{hu2020}
J.~Hu, X.~Liu, Z.-W. Wen, and Y.-X. Yuan.
\newblock A brief introduction to manifold optimization.
\newblock {\em Journal of the Operations Research Society of China},
  8(2):199--248, 2020.

\bibitem{jiang2019}
R.~Jiang and D.~Li.
\newblock Second order cone constrained convex relaxations for nonconvex
  quadratically constrained quadratic programming.
\newblock {\em J. Global Optim.}, 75(2):461--494, 2019.

\bibitem{kim2020}
J.~Kim, M.~Kang, D.~Kim, S.-Y. Ha, and I.~Yang.
\newblock A stochastic consensus method for nonconvex optimization on the
  stiefel manifold.
\newblock In {\em 2020 59th IEEE Conference on Decision and Control (CDC)},
  pages 1050--1057. IEEE, 2020.

\bibitem{More.Sorensen.1983}
J.~J. Mor{\'e} and D.~C. Sorensen.
\newblock Computing a trust region step.
\newblock {\em SIAM J. Sci. Statist. Comput.}, 4(3):553--572, 1983.

\bibitem{nemirovski2007}
A.~Nemirovski.
\newblock Sums of random symmetric matrices and quadratic optimization under
  orthogonality constraints.
\newblock {\em Math. Program.}, 109(2-3, Ser. B):283--317, 2007.

\bibitem{overton1992}
M.~L. Overton and R.~S. Womersley.
\newblock On the sum of the largest eigenvalues of a symmetric matrix.
\newblock {\em SIAM Journal on Matrix Analysis and Applications}, 13(1):41--45,
  1992.

\bibitem{smith1993}
S.~T. Smith.
\newblock {\em Geometric optimization methods for adaptive filtering}.
\newblock PhD thesis, Harvard University, Cambridge, MA, May 1993.

\bibitem{wen2013}
Z.~Wen and W.~Yin.
\newblock A feasible method for optimization with orthogonality constraints.
\newblock {\em Mathematical Programming}, 142(1):397--434, 2013.

\bibitem{zhang2014}
X.~Zhang, J.~Zhu, Z.~Wen, and A.~Zhou.
\newblock Gradient type optimization methods for electronic structure
  calculations.
\newblock {\em SIAM Journal on Scientific Computing}, 36(3):C265--C289, 2014.

\end{thebibliography}

\end{document}